\documentclass[a4paper,10pt]{article}
\usepackage{hyperref}

\usepackage[utf8]{inputenc}

\usepackage[english ]{babel}
 
\usepackage{xcolor}
\usepackage{bbding}

\usepackage[utf8]{inputenc}

\usepackage{tikz-cd}

\usepackage{amsfonts}
\usepackage{amssymb}
\usepackage{graphicx}
\usepackage{mathtools}
\usepackage{skak} 
\usepackage{foekfont}

\usepackage[T1]{fontenc}

\usepackage{amsmath,mathdots}

\usepackage[utf8]{inputenc}
\usepackage{devanagari}

\usepackage{mathtools}
\usepackage{mathdots}
\usepackage{tikz}

\definecolor{amber(sae/ece)}{rgb}{1.0, 0.49, 0.0}
\newfont{\rsfsten}{rsfs10 scaled 1200}

\usepackage{MnSymbol}
\usepackage{tikz}
\usepackage{amscd}
\usepackage{MnSymbol}
\usepackage{wasysym}

\usepackage{mathrsfs}

\newcommand*{\rom}[1]{\expandafter\@slowromancap\romannumeral #1@}

\usepackage{xcolor}

\usepackage{wrapfig}
\newcommand{\tightunderset}[2]{%
  \mathop{#2}\limits_{\vbox to .3ex{\kern-0.95ex\hbox{$#1$}\vss}}}

\usepackage[utf8]{inputenc}

\usepackage{CJKutf8}

\usepackage{pmboxdraw}

\usepackage{stmaryrd}
\usepackage{url}
\usepackage{harmony}
\usepackage{ifsym}

\usepackage{eufrak}

\usepackage{rotating}

\usepackage{dingbat}

\makeatother

\title{ Convex Polytopes, Dihedral Angles, Mean Curvature  and Scalar Curvature}
\author{Misha Gromov, }

\begin{document}
  \maketitle

\begin{abstract} 

We approximate boundaries of convex  polytopes $X\subset \mathbb R^n$ by smooth hypersurfaces $Y=Y_\varepsilon$ 
with  {\it positive mean curvatures} and,  by using  basic geometric relations  between 
the scalar curvatures of Riemannian manifolds  and the mean curvatures of their boundaries,
establish   {\it  lower  bound  on   the dihedral angles} of     $X$. 
 \end{abstract} 
\tableofcontents

\section  
{Combinatorial Spread, $\square_\rangle$-Spread and  $\square_\rangle$-Inequality}

 Let  $X\subset \mathbb R^n$ be  a compact  convex polytope let $\partial X$ denote its (topologically spherical) boundary and let  $X^{\circledcirc}$ be {\it the dual convex   tessellation} of the sphere 
$S^{n-1}$, i.e.  where  $(n-k-1)$-cells are the sets  of the (unit normal to the)  supporting hyperplanes to $X$ along the interiors of the $k$-faces   of $X$. 

Let $E= E(X)\subset S^{n-1}$ be the edge  graph of $X^{\circledcirc}$.   Combinatorially, this is 
 the {\it $(n-2)$-adjacency  graph},  where 
 the set  of the  {\it $(n-1)$-faces} $F$ of $X$  is taken for the set  vertices and 
 where  the edges $e$ in  $E$ correspond to the pairs of  $(n-2)$-{\it adjacent faces}:

  {\sf vertices $v_1$ and $v_2$  are joined by an edge  $e= e_{12}$, whenever the corresponding {\it closed faces} $\bar F_{1}, \bar F_{2}\subset X$  meet over   a  closed  
 $(n-2)$-face, namely $\bar F_{12}= F_{1}\cap F_{2}\subset X$.}

{\it Remark on Adjacency and on  Simple  Polytopes.} Recall that a  convex polytope  $X$ is {\it simple} if  \vspace {1mm}

 \hspace {33mm} adjacent $\implies$ $(n-2)$-adjacent,  \vspace {1mm}

 \hspace {-6mm}   where "adjacent" signifies that the intersection $F_1\cap F_2$ is {\it non-empty}, i.e. $F_1$ and $F_2$ meet at a 
 vertex in $X$.

The {\it combinatorial distance} $dist_{comb} (F_1,F_2)$ is the length of the shortest path in $E$ between the corresponding vertices corresponding to  $F_1$ and $F_2$. 

For instance, these distances   between opposite faces in the   n-cube $[-1,1]^n$ are equal to 2. 
\vspace {1mm}

Let $\angle_{1,2}=\angle( F_1,F_2)$ denote the dihedral angle between  $(n-2)$-{\it adjacent faces}  and let $\rangle$ stands for the complementary 
angle, 
$$\rangle_{1,2}=\pi-\angle_{1,2},$$ 
that     is  the spherical  arc length of the edge  $e_{12}\subset  S^{n-1}$  dual to the $(n-2)$-face 
$F_{12}=F_1\cap F_2$.

{\it  \textbf{ $\rangle$-Angular  Distance.}}  The {\it  angular distance}
 or {\it $\rangle$-distance}  $dist_\rangle(F_1,F_2)$  between (not necessarily $(n-2)$-adjacent)  $(n-1)$-faces $F_1$ and $F_2$ 
  in $X$ is the minimum of the  spherical  lengths of  edge paths in $E$ between the  vertices of $X^{\circledcirc}$  dual to these faces.

Accordingly, the {\it $\rangle$- (angular)  distance} between (unions of) sets of faces, say $\mathcal F_1, \mathcal F_2\subset V$, is the minimum 
of the $\rangle$-distances between the faces in these sets, 
$$dist_\rangle(\mathcal F_1, \mathcal F_2)=\min_{F_1\in \mathcal F_1, F_2\in\mathcal F_2}dist_\rangle(F_1,F_2).$$

{\it Cubical Example.} The $\rangle$-distances between  opposite faces of the $n$-cube $\square^n=[-1,1]^n$ are equal to $\pi.$

{\it Combinatorial and Angular Spreads}. Let $\square^k_{comb}(X)$ be the {\it maximum of the numbers }   $d\geq 0$, such that 
$X$ admits a continuous map to the  $k$-cube,  
$$\Phi:X\to \square^k=[-1,1]^k$$ 
with the following properties.\footnote {If no such map exists, then 
let $\square^k_{comb}(X)=0$.}

$\bullet_{comb}$  The $\Phi$-pullbacks of the  $(k-1)$-faces from  $\square^k $ are unions of $(n-1)$-faces in  $X$.

$\bullet_{dist}$  The combinatorial distances between the pullbacks of the opposite  cubical faces $\underline F_{i\mp}\subset \square^k$  are $\geq d$,
$$dist_{comb}(\Phi  (\underline F_{i-} ), \Phi  (\underline F_{i+}))\geq  d,  \mbox { }  i=1,...,k. $$

$\bullet_{deg}$ The the induced relative  homology homomorphism 
$$ \Phi_\ast:H_k(X, \Phi^{-1}(\partial\square^k))\to H_k(\square^k,\partial\square^k)=\mathbb Z$$
{\it doesn't vanish.}

(If $k=n$, this is equivalent to    $\Phi^{-1}(\partial\square^n)=
\partial X$ and to non-vanishing of the degree of the  map $\Phi:\partial  X\to\partial\square^n$. For instance, homeomorphusms
$\Phi:X\to\square^n$ satisfy tis condition.\footnote{The topological  degree is defined for all {\it continuous  equidimensional}  maps $f$ between {\it oriented manifolds}, e.g. such as our spherical $\partial X$ and $\partial \square^n$, where the non-vanishing condition $deg(f)\neq 0$  doesn't depend on the orientation for connected (orientable) manifolds.
Also the   degree is defined for  the boundary respecting maps between manifolds with boundaries. }
)

\vspace{1mm} 

Similarly define the angular spread $\square^k_{\rangle}(X)$ with the $dist_\rangle$ inequality instead of   $dist_{comb}$: 
$$dist_{\rangle}(\Phi  (\underline F_{i-} ), \Phi  (\underline F_{i+}))\geq  d,  \mbox { }  i=1,...,k. \leqno {\hspace {6mm} \bullet_{dist_\rangle}}$$

Observe that the combinatorial and   the angular spreads satisfy
$$  \square^{n}(X)\leq  \square^{n-1}(X)\leq...  \leq  \square^{1}(X)=diam(X),$$
where the diameter refers to the combinatorial and to the  angular distances correspondingly, 

and that 
$$  \square^{n}(X)\geq{\square^{n-1}(X) \over 2n+2}\mbox{ for    all convex $n$-polytopes $X$}.$$

{\it \textbf {$\square^3_N$-Example.}} Let $\square ^3_N$ be the subdivision of  the  3-cube 
$\square ^3=\square^3_1$, where each  2-face is subdivided into 
$N^2$  equal squares  in an obvious way.  (If you wish it to became simple,  $\varepsilon$-perturb  with $\varepsilon << 1/N$  the  edges of) these  small squares,   such 
 that the resulting  subdivision $\square^3_{N, \varepsilon}$ has   three   squares  at each vertex.)
 Then  the combinatorial  $\square^3$-spread of the so subdivided cube   is   $N+1$.

{\it "Random" Example.} Apparently, the  combinatorial $\square$-spread of  a suitably defined  {\it random $n$-polytope} with $M$ faces 
  (see section   8)    
       grows, roughly,  as  $\sqrt [n-1] {M}$.

\textbf {1.A. Angular Spread  Theorem.}  {\it The top-dimensioanl   $\rangle$-spreads, 
 of 
all  compact convex   $n$-polytopes $X\subset \mathbb R^n$ are bounded by a universal constant,
$$\square^n_\rangle(X)\leq D=D_n\leq  2(n-1)\sqrt {n}.$$}

 We shall proof this in section 5 by reduction to the  {\it  normalized mean curvature  mapping theorem} 
 (see section 2)  the proof of which    (see section  2.1)      depends on    the index theory for Dirac operators on Riemannian spin manifolds with positive  scalar curvatures (see sections 3.1.2 and  3.5 in [Gr2021]).

{\it \textbf { 1.B. Corollary.}} The minimum of the complementary angles of $X$ is bounded by the  combinatorial spread  
$\square^n _{comb}(X)$
as follows,
$$\rangle_{min}(X)\leq D\frac{1}{\square^n_{comb}(X)},$$
for the  above constant $D$.

{\it \textbf { 1.C. Conjecture.}} {\sf The above $D=D_n$  is equal to $\pi$.\footnote{In view  of what Karim Adiprasito recently told me, this  $D$  is better to be the one from 1.C  rather  than from 1.B. }}

  \vspace{1mm}

{\it Remark.} Probbaly 1.A, 1.B and 1.C generalize to
all convex tessellation  of $S^{n-1}.$\footnote {There must exist convex tessellations, which  don't  come  from convex polytopes, 
but I failed to find a reference or  to find such tessellations  by myself.}  (See  next  section for more about it.)

{\it \textbf {Acknowledgements.}} I am grateful to J\'anos Pach for his interest in this article  and to  Christina Sormani and Karim Adiprasito  who pointed out several errors in an earlier version of the paper.


\subsection  {Combinatorially Large Polytopes with  Large Complementary Angles   }  

 Dirac  operators notwithstanding, evaluation of the  ranges of 
   possible  values of the dihedral angles of 
    polytopes depending on their  combinatorial types and/or determination of the combinatorial and metric 
  geometries of polytopes with  all  complementary  angles   $\rangle(X)$  bounded from below 
  remains problematic.

It is known here (Steinitz?) that if 
 $\rangle_{min}(X) \geq \pi/2$,  then $X$ is  the  product of simplices.
But -- this was 
pointed out to me by Karim Adiprasito  three years ago --  there is no bound on the number of
 faces of $X$ for $\rangle_{min}(X)\geq \alpha $ for small $\alpha>0$. Later,  I found the following on the web.

{\it ${\pi}/4$-Example.} Chop off the corners 
  from the prism $\octagon_i\times [0,\delta]\subset \mathbb R^3$, where $\octagon_N$ is the regular $N$-gon, $N=3,4,.... $  and $\delta=10^{-N}$, such that 
 this "chopping"   fully consumes the $\delta$-edges of the prisms, and such that all  complementary  dihedral angles of the resulting polytopes $X_N$ are mutually equal and  satisfy $\alpha_N\to \pi/4$  for $N\to \infty.$
 see 
 [MSE]. 

 Recently,  Karim informed  me [AK]  that infinity of  combinatorial types of convex polytopes with $\rangle_{min}(X)\geq \alpha $ exits if an only if      $\alpha<\pi/3$. 
 
 Below, in a similar spirit, we   construct  $n$-polytopes,  which

 {\it  "infinitely stretch" in $n-2$ directions, while having all complementary dihedral angles  bounded from below.}

{\it \textbf {1.1.A.  Skyscrapers}}. Given convex polytopes $\mathbf 0 \in X_1\subset ...\subset X_N\subset  \mathbb R^{k}$
and numbers  $  h_1>...>h_N\geq 0$ let
 $$\bigdoublewedge^N\{X_i\}=\bigdoublewedge^N_{\{h_i\}}\{X_i\}\subset   \mathbb R^{k}\times \mathbb R_+\subset  \mathbb R^{k+1}$$
denote the  intersection of the cones of heights $h_i$ over $X_i$,  where  the top vertices of these cones lie on the "vertical" axes
$\mathbf 0\times \mathbb R_+\subset  \mathbb R^{k}\times \mathbb R_+$,
$$\bigdoublewedge^N\{X_i\}=\bigcap_i^N cone_{h_{i}} (X_i).$$

Such a $\bigdoublewedge=\bigdoublewedge^N=\bigdoublewedge\{X_i\}$ is called {\it a skyscraper} with the {\it bottom} $X_1$ and the {\it  top} $X_N$ if the following holds:

 $\bullet_\cap$ the bounary of $\bigdoublewedge\{X_i\}$ has {\it non-empty intersections}  with {\it all}  open  $(n-1)$-faces of the cones $cone_{h_i} ( X_i^{h_i})$ and

   {\it two closed side  faces  $F_1,F_2\subset\bigdoublewedge\{X_i\} $ {\it do not} intersect unless they contained in faces of the cone
over some $X_i$,
 $$F_1,F_2\subset \partial cone_{h_i} ( X_i),$$  
 or in faces of  two consecutive cones.}
$$F_1 \subset \partial cone_{h_i}  (X_i)\mbox {  and } F_2 \subset \partial cone_{h_{i+1}} (X_{i+1}).$$

 Notice that $\bullet_\cap$ implies the following:

$\bullet_\#$ the number of the $k$-faces of $\bigdoublewedge\{X_i\}$ satisfies 
$$\#_{k}\bigdoublewedge\{X_i\}=1+\sum_{i=1}^N\#_{k-1}X_i.$$

  Observe that $\bullet_\cap$ can be  achieved with suitable $h_i$   and homothetically scaled $X_i$.
 
$\bullet _{\{\lambda_i\}}$ {\it Given $X_i$  \  there exist $h_1>...>h_i>...>h_N$ and $0<\lambda_1< ...<\lambda_i<...<\lambda_N,$
such that $\bullet_\cap$ is satisfied by 
$\bigdoublewedge^N_{\{h_i\}}\{\lambda_iX_i\}$.\footnote {$\lambda X=\{\lambda x \}_{x\in X}\subset \mathbb R^k$.} }

The usefulness of his for our purpose is due to the following obvious property  of skyscrapers.

{\it \textbf {1.1.B.  Large $\rangle $   Lemma}.} {\sf  Let $\bigdoublewedge_{\{h_i\}}\{X_i\}$ be  a  {\it a skyscraper} (with the  bottom $X_1$ and the top $X_N$), such that
 the complementary angles of all $X_i$ as well as 
 (by definition acute)  angles  between the pairs of
 hyperplanes, which  define the    faces of consequtive    $\underline X_i$  and of  $\underline X_{i+1}$,
 are strictly   bounded from below by $\alpha>0$.

Then there exist  a (large)  positive  number $C$ such that vertically $C$-stretced $\bigdoublewedge$, 
that is 
$C\uparrow\bigdoublewedge=\bigdoublewedge_{\{C\cdot h_i\}}\{X_i\}$ }  (a true skyscraper) {\sf 
has the  complementary dihedral angles between the side faces bounded from below by $\alpha$  
 $$\rangle_{side} \left(C\uparrow\bigdoublewedge\right) >\min (\alpha, \pi/2),$$ }
while these angles  at the bottom face  are  $>\pi/2$.
\vspace{1mm}

{\it $\pi/3$-Example.} Let  $ X_i\subset \mathbb R^2$  be regular triangles,
 where $ X_i=\underline X_{i+2}$  and $X_2=-X_1$. 
 Then  the complementary side dihedral angles $\rangle $ of the  corresponding  skyscraper   
 $$X_{\hexagon}=X_{\hexagon}(N,C)=C\uparrow\bigdoublewedge=\bigdoublewedge_{\{C\cdot (N-i+1)\}} \{(N+1)^iX_i\}$$
 satisfy 
 $$\rangle_{side} (X_{\hexagon}(C,N))\to \pi/3\mbox {  for } C\to \infty,$$
  
  while  $$ \rangle_{bottom} (X_{\hexagon}(C,N))\to \pi/2.$$
 
 (This,  I guess,  must be exactly  Adiprasito's example.) \vspace{1mm}

{\it Remarks.} (a) 
 The  essential 
 difference of $X_{\hexagon}(C,N))$  from the above "pruned"   prism  $\octagon_N\times [0,\delta]$ is that
 
 \hspace{1mm}  {\it  the combinatorial diameters of $X_{\hexagon}(C,N)$ tend to infinity  for $N\to\infty$.}
 
  In fact, \vspace{1mm}

\hspace {40mm} $diam_{comb}( X_{\hexagon}(C,N))=N+1.$

(b)    The Cartesian products of $m$ copies  of  $X_{\hexagon}$  provide  examples  of 
$3m$-polytopes,
 $m=1, 2....$, with all complementary angles $\geq \frac {\pi}{3}-\varepsilon$ for all $ \varepsilon>0$
and with arbitrarily large $\square_{comb}^m$-spreads.

(c)  {\it The directional limit  set of the faces} of the 3-polytops $X_{\hexagon}(C,N))$  for $N,C\to \infty$, that is the Hausdorff limit  of the  sets of vertices of the  dual tessellations  $X^{\circledcirc}_{\hexagon}(C,N))$
 of $S^2$,  is a 7-point set:  a regular hexagon on the equator plus the south pole, while 
 similar $X_N$  with suitably  rotated  triangles $\underline X_i$ may have arbitrary limit  sets on the equator.

{\it Question.} Is this limit set  always {\it discrete} away from an equatorial circle $S^1\subset S^2$?

(Adiprasito bound  $ \rangle_{min}\leq\pi/3$ makes this plausible     for 
$ \rangle_{min}\underset {N\to \infty}\to\pi/3$.)
\vspace {1mm}
\vspace{1mm}

{\it \textbf {1.1.C.  Skyscrapers on Skyscrapers}.} Finiteness of the directional   limit sets of Skyscraper $\bigdoublewedge\{X_i\}$ allows a lower bound on the complementary dihedral angles of  
 double skyscraper  $\bigdoublewedge\{ \bigdoublewedge_j\{X_i\}\}$, etc.

{\it $\pi/3(2n-5)$-Example.} Let $\rho_{n,k}(\triangle)\subset \mathbb R^2$, $n=3,4,...$, $k=0,...2n-5,$ be the regular triangle rotated by $\rho_{n,k}=k\pi/3(2n-5)$ and define by induction on $m$ polytopes 
$\mathcal X_m=\mathcal X_m(n, N_{n,m-2}, C_{n, m-2})\subset \mathbb R^{m}=\mathbb R^2\times \mathbb R_+^{m-2}$, $m=3.4,...n$, as follows.  

Let 
$$\mathcal X_3=\bigdoublewedge \{\rho_{n,0}(\triangle),  \rho_{n,1}(\triangle)|N_{n,1}\}=
C_{n,1}\uparrow\bigdoublewedge_{\{h_i\}}^{2N{n,1}} \{\lambda_i| \rho_{n,0}(\triangle),   \rho_{n,1}
(\triangle)|N_{n,1}\},$$
where
$\{\lambda_i|A,  B|N\}$ stands for $\{\lambda_1A, \lambda_{2}B,...\lambda_{2N-1}A, \lambda_{2N}B\}$
and where the constants $h_i$ and $\lambda_i$ are chosen as in 1.1.A and where eventually  $C_{n,1}\to \infty$ as earlier.

Then we slightly modify $\mathcal X_3$ by turning the base  2-face $F_{base}=\triangle=\mathcal X_3\cap \mathbb R^2\times 0$    by $\pi/4$, call the result  $\mathcal X'_3$ and inductively define 
$$ \mathcal X'_{m+1} =\bigdoublewedge \{  \rho_{n,2(m-3)} (\mathcal X'_m),  \rho_{n,2(m-3)+1} (\mathcal X'_m)|N_{n,m+1}\},$$
where the rotations $\rho$  apply to the $\mathbb R^2$-factor in 
$\mathbb R^2\times \mathbb R_+^{m-2}\supset \mathcal X'_m$
and where the implicit $h,\lambda$ and $C$- constants are adjusted as earlier.    

It is easy to show -- we leave checking this to the reader that

{\sf {\huge $\star$}\hspace{-0.5mm}$\square$\hspace{1mm}  {\it the $\square^{n-2}_{comb}$ stretch} of   $\mathcal X'_n$ can be made  {\it arbitrarily  large}   with all $N_{n,m}\to \infty$\footnote {The  boundary of $\mathcal X'_n$ (as well as that of $\mathcal X_n$)  with the $dist_{comb}$-geometry is shaped  roughly the same as
({\it properly coarsely homotopy equivalent to})  the rectangular solid $[0,N_{n,1}]\times...\times
  [0,N_{n, n-2}]$.  }   

and that 

{\huge $\star$}\hspace{-0.5mm}$\rangle$\hspace{1mm} {\it the  complementary  dihedral angles }  of  $\mathcal X'_n$ satisfy
$$\rangle( \mathcal X'_n)\geq \frac{ \pi}{3(2n-5)}-\varepsilon,$$
where $\varepsilon>0$ can be made {\it arbitrarily small} with  $C_{n,m}\to \infty.$ }

\vspace {1mm}

Probbaly  a  skyscraper pattern  is present in all polytopes $X$
with    $\square^{n-2} _{comb}>> \frac{1}{\rangle_{min}}$.  We partly  justify this (conclusively only for  $dim(X)= n=3$) 
by looking 
at the the dual spherical tessellations $X^{\circledcirc}$ as  follows.
 
 Given a cellular   tessellation $T^{\circledcirc}$, e.g. a triangulation,  of an   $(n-1)$-manifold $Y$, define the combinatorial distance between cells, as earlier, by the lengths of  
minimal chains of cells,  denote  this by $dist_{\circledcirc}$ and  define the combinatorial 
$\square^k_{\circledcirc}(Y)=\square^k_{comb}(T^\circledcirc)$,  including $diam_{\circledcirc}=\square^1_{\circledcirc}$, via continuous maps $\Phi:Y\to \square^k$   by just saying "cell" instead of "face".

{\it \textbf {1.1.D. Large Subdomain  Lemma.}} {\sf Let 
$T^\circledcirc$  be a convex tessellation of the unit sphere $S^{n-1}$, where the cells are called $\triangle$,  and let  $B^{\circledcirc}_s\subset S^{n-1}$, $s\in S^{n-1}$, denote the union of closed cells which contain $s$.

Then, given a (small) number  $v>0$,} 

{\sl there exists a connected cellular (i.e. a union of cells)   domain  $U^{\circledcirc}$ in the sphere $S^{n-1}$, 
 such that the spherical volumes of the "${\circledcirc}$-balls" $B^{\circledcirc}_s$ around all points in $U$ are bounded by  
$$vol(B^{\circledcirc}_{s})\leq v, s\in U^{\circledcirc},$$
and such that the $\square_{comb}$-spreads of  $U^{\circledcirc}$  are bounded from  below  by these of $T^\circledcirc$ as follows:
$$\square^k_{\circledcirc}(U^{\circledcirc})\geq {const }\cdot{v} \cdot\square^k_{\circledcirc}(U^{\circledcirc})-1,
 \mbox { }const=const_n>(10 n)^{-n}.$$}

Indeed, the cardinalities   $N=N(v)  $  of   subsets $S\subset S^{n-1}$ of "$v$-{\it thick}" points  $s\in S^{n-1}$, i.e. with
$vol(B^\circledcirc_s)\geq v$,  such that  no pair of  different points from $S$ is contained  in the same closed  cell of $T^{\circledcirc}$,   are bounded by $N=\frac {vol(S^{n-1})}{v}$, while the combinatorial diameters of all $\circledcirc$-balls, are at most  2, 
$$diam_{comb}(B^\circledcirc_s) \leq 2\mbox { for al $s\in S^{n-1}$}.$$

Therefore, given a map $\Phi =\{\Phi_1,...,\Phi_k\}: S^{n-1}\to [-1,1 ]^k$ from the definition of $\square_{comb}^k$, there exist  gaps between
  pairs 
 of  {\it  neighbouring mages},  say   $t_i=\Phi_i(s_i) , t_i'=\Phi_i(s'_i) \in [-1,1]$,  $i=1,...k$, of  pairs of "$v$-thick"  vertices $s_i$ and $s_i'$, 
 such that 
$$dist_{comb} (\Phi_i^{-1}[-1,t_i], \Phi_i^{-1}[t_i',1])\geq \frac {1}{N}dist_{comb}(\Phi_i^{-1}(-1), \Phi_i^{-1}(1))-2(N+1)$$
and the "$B^{\circledcirc}$-enlargement" of the intersection $U_\cap$ of the  pullbacks   $\Phi^{-1}[t_i,t_i']\subset S^{n-1}$ is  taken for the required $U^{\circledcirc}$
$$U_\cap=\bigcap_{i=1}^k\Phi^{-1}[t_i,t_i']  \mbox{   and }  U^{\circledcirc}=
\bigcup_{s\in U_\cap}B^{\circledcirc}_s.$$

Here is another   obvious  observation.

{\it \textbf {1.1.E. Narrow Band  Lemma.}}
{\sf If  the edges (1-cells) from  $T^{\circledcirc}$, adjacent to a vertex $s\in S^{n-1}$, 
have  lengths $\geq l$ and if $vol(B^{\circledcirc}_s)\leq v$,
then $B_s^{\circledcirc}\subset S^{n-1}$ is {\it contained in the $\delta$-neighbourhood of an equatorial
sphere} $S^{n-2}\subset S^{n-1}$, 
$$B^{\circledcirc}_s\subset U_\delta (S^{n-2}),$$
where this $\delta=\delta_n(l,v)>0$  satisfies  for all $n$ and $l>0$,  
$$\delta_n(l,v)\to 0 \mbox { for }v\to 0.$$
Moreover, if all $(n-2)$-cells $\triangle^{n-2}$ adjacent to $s$ have
$$vol_{n-2} (\triangle^{n-2})\geq a>0,$$
then this equatorial $S^{n-2}\subset S^{n-1}$ is {\it unique up to an $\varepsilon$-perturbation}, i.e. all 
equators for which $ U_\delta (S^{n-2}) \supset B^{\circledcirc}$ lie within distance   $\varepsilon$ one from another, where
$$\varepsilon=\varepsilon_n(v,\delta, a)\to 0\mbox {  for }v\to 0.$$}

{\it \textbf {1.1.F. Corollary: Elementary Bound on $\square^2_{comb}$.}} {\sf  A  lower bound by $a>0$ on the $(n-2)$-volumes of  $(n-2)$-cells in  a convex tessellation  $T^{\circledcirc}$ of $S^{n-1}$  implies
 an upper bound on the combinatorial $\square^2$-spread  of  $T^{\circledcirc}$,
$$ vol_{n-2}(F^{n-2})\geq a>0\implies \square^2  (T^{\circledcirc})\leq\Theta_n(a).$$
 where $\Theta_n$   
is a (bounded monotone decreasing)   function in  $a>0$.}

(If $n=2$, this is just  a qualitative version of  1.B.)

{\it  Proof}.  It follows from   1.1.F. and  1.1.E that  the above
{\sf $U^{\circledcirc}\subset S^{n-1}$ is contained in a $\delta'$-neighbourhood  of an equator  
$S^{n-2}\subset S^{n-1}$, where, for a fixed $a>0$, 
$$\square^k_{comb}(T^{\circledcirc})\to \infty\implies  \square^k_{comb}(U^{\circledcirc})\to \infty
  \mbox {  for all   $k=1.2...$  .}$$}

Then one sees that,  for $\delta'$ much smaller than $a$, this  $U^{\circledcirc}$ admits a  cellular map of {\it degree one}  from the cylinder 
$S^{n-2}\times [0,1]$, which is   decomposed into $m\times M$ cells,   which are products of cells of some triangulation of  $S^{n-2}$ into $m$-simplices and a  decompositions of $[0,1]$
 into $M$ segments,  where $m$ is bounded by  
a constant depending on $a$.

 It follows that  $\square^2( U^{\circledcirc})$  is also $ \leq m$, hence, it is  bounded in
  terms of $a>0$. QED.

{\it  Remark.} The above shows that if   $\square^1_{comb}(T^{\circledcirc})\to \infty$ with the $(n-2)$ volumes of all  $(n-2)$-cells bounded from below by $a$, then the unit sphere $S^{n-2}$ acquires several limit  tessellations with the same bound on the volumes of their $(n-2)$-cells  and some cells spanned by vertices of different tessellations.

Then, for instance, by looking on  pairs of such tessellations,  one   recovers    a special case of   Adiprasito's result  for $n=3$.

{\it \textbf {1.1.G.   Conjecture.}}  {\sf  For all $k=1,... n-2$,  a  lower bound   on the $k$-volumes of the $k$-cells in  a convex tessellation  $T^{\circledcirc}$ of $S^{n-1}$  implies
 an upper bound on the combinatorial $\square^{n-k}$-spread  of  $T^{\circledcirc}$.}
 
 Conversely,  
 
 {\sf there exist convex tessellations  $T^{\circledcirc}$ of $S^{n-1}$ with arbitrary large  $\square_{comb}^{n-k-1} ( T^{\circledcirc})$ and  with the volumes of all $k$-cells bounded away from zero.

Moreover, there are such  $T^{\circledcirc}$, which are dual of convex polytopes  $X\subset \mathbb R^n$.}
 (A quantitative form  of a special case  of this  is suggested in 6.B.)


\section{ Manifolds with Corners, Mean Convexity and Distance $dist_\rangle^\natural$.}

Let  $X$ be  a smooth $n$-manifold with corners,  i.e.   locally,
   at all $x\in X$, it is  diffeomorphic to a convex polytope  $Q=Q_x\subset \mathbb R^n$. 

For instance,  diffeomorphic images of convex polytopes are manifolds with corners.

Also recall that the {\it mean curvature} of a cooriented hypersurface in  a Riemannian manifold is the sum
of the principal  curvatures.

{\it Example.}  The $R$-sphere $S^{n-1}(R)\subset \mathbb R^n $  and the round cylinder   
 $S^{n-2}\times \mathbb R^1\subset \mathbb R^n $ satisfy 
$$mean.curv(S^n)=\frac{n-1}{R} \mbox { and  }  mean.curv (S^{n-2}(R)\times \mathbb R^1)= \frac{n-2}{R}.$$

A Riemannian manifold with corners is called {\it mean convex} if all its $(n-1)$
-faces $F\subset \partial X$ have non-negative 
mean curvatures.

For instance, convex domains in $\mathbb R^n$ with corners are mean convex.

Given a smooth curve in the boundary of a manifold with corners, say   $\gamma \subset \partial X$,   which  doesn't intersect $(n-2)$-faces of $X$ and 
which meets   all $(n-2)$-faces of $X$ transversally, 
say at  $x_i\in \partial X$, $i=1, ...,  j$,  let 
$$length_\rangle^\natural(\gamma) = length^\natural(\gamma) + \sum_{i=1}^j \rangle_{x_i},$$
where $length^\natural(\gamma)=\int_\gamma mean.curv(\partial X))d\gamma$, 
where  $ \rangle_{x_i}$ are the complementary dihedral   angles, $\rangle_{x_i}=\pi-\angle_{x_i}$ and  where   the  dihedral angle $\angle_{x_i}$  of $X$ at the point $x_i$ is  the angle between the (naturally cooriented) 
$(n-1)$-dimensional tangent  spaces $T_i, T'_i\subset T_{x_i}(X)$ to the two $(n-1)$-faces   adjacent to the (n-2)-face, which contains  $ x_i$.

Next, assuming  $X$  is mean convex and $x_1,x_2\in\partial X$ are  contained {\it inside} $(n-1)$-faces,  let 
   $$ dist_\rangle^\natural(x_1, x_2)=\inf_{\gamma_{1,2}} length_\rangle^\natural(\gamma_{1,2}), $$
 where the infimum is taken over  the above kind of curves 
$\gamma_{1,2}\subset \partial X$ between $x_1$ and  $x_2$. 

Although this $dist_\rangle^\natural$ is defined not for all points and it   may  vanish at some pairs  of non-equal points,  we treat it as a true distance; in particular,  we define 
the corresponding distance between $(n-1)$-faces\footnote {An open   $k$-face in $X$ is understood as a maximal  {\it connected} subset  where $X$ is 
locally diffeomorphic to a polyhedral $k$-face.}  $F_1,F_2\subset \partial X$ in the usual way:
$$ dist_\rangle^\natural(F_1, F_2)=\inf_{x_1,x_2} dist_\rangle^\natural(x_1,x_2)\mbox
  {  for $x_1\in F_1$  and $x_2\in F_2$.} $$

If the  boundary of $X$ contains no corners, i.e. it is  smooth,  then the corresponding distance is denoted $dist^\natural$. This is a true positive distance if $X$ is {\it strictly mean convex}, i.e.  $mean.curv(\partial X)>0.$ 
\vspace{1mm}

{\it Semi(in)stability of $dist^\natural$}. {\sf An arbitrarily  $C^1$-small perturbation of a smooth convex hypersurface 
$ Y\subset \mathbb R^n $, $n\geq 3, $    may significantly diminish the metric  $dist^\natural$  on $Y$.} 

For instance, 

{\sl  the unit sphere $S^{n-1}\subset \mathbb R^n$, which has $diam^\natural(Y)=\frac {(n-1)\pi}{2}$,   can be $C^1$-approximated by smooth convex hypersurfaces 
$Y_\varepsilon$ with $diam^\natural(Y_\varepsilon)= {\pi}+ \varepsilon$  for all $ \varepsilon>0$ as follows.}\footnote {This can't be done with $\varepsilon=0$.} 

Let $A=A_{N,\delta} \subset S^{n-1} $, where $N>\delta ^{-2n}$, be the union of regular equatorial $N$-gons in general position,  such  
for all pairs of points in $S^{n-1}$, there are our $N$-gons  passing $\varepsilon$  close to both points.    Let   $B_A\subset \mathbb  R^n$ be the intersection of subspaces bounded by the  hyperplanes tangent to $S^{n-1} $ at the vertices of the $N$-gons and  let $Y(N, \delta, \epsilon) $ be the boundary of the 
$\epsilon $-neighbourhood   of $B_A$  for $0<\epsilon \leq N^{-2n}$.

Then  $$diam^\natural(Y(N, \delta, \epsilon))\to \pi \mbox  { for }\delta \to 0,$$

and, for the same  reason, 

{\it all  $\natural$-spreads of $Y(N, \delta, \epsilon)$  converge to the ordinary speads of the unit sphere, 
$$  \square^k_\natural (Y(N, \delta, \epsilon)\to  \square^k (S^{n-1})=
 \frac {1}{n-1}  \square^k_\natural (S^{n-1}),$$}
 where, observe for instance,    $\square^n (S^{n-1})=2\arcsin {1\over \sqrt n}.$\footnote {Here $\square^n (S^{n-1})$ means  $\square^n (S_+^{n})$, i.e. this is defined via maps $S^{n-1} \to \partial \square^n$ of positive degrees.}

\vspace {1mm}

But the metric $dist^\natural$ of  a compact convex hypersurface $Y$  can't {\it everywhere}   significantly increase 
under  small 
$C^0$-perturbations of $Y$.

In fact, if $Y=S^{n-1} $, this  follows from theorem  2.A  below, which is a special case of  the {\it normalized
mean curvature mapping theorem}  from section 3.5 of [Gr2021] and which makes the key
 ingredient  of the proof of 1.A.

\textbf {2.A. Euclidean  $dist^\natural$-Non-Contraction Theorem.} {\sf Let $X$ be a compact oriented {\it mean convex} Riemannian $n$-manifold with smooth boundary, let $B\subset \mathbb R^n$ be a smooth compact  {\it convex} domain. e.g the unit ball,  and let 
$f:\partial X\to \partial B$ be a smooth map, which 
  which has  {\it non-zero degree}. 

If $X$  has {\it non-negative} scalar curvature, $Sc(X)\geq 0$,  and if  $X$ is {\it spin},\footnote {{An oriented vector bundle  $T \to X$ is {\it  spin} if the associated principle bundle $G\to X$ with the fibres $G_x=SL(n, \mathbb R)$,  $n=rank(T)$,  admits a double covering $s: \tilde G\to G$, such that  the pullbacks of the fibers $s^{-1}(G_x)$ are {\it connected};  
   an orientable  manifold $X$ is {\it spin} if the tangent bundle $T(X)$   is spin.
  
  A necessary and sufficient condition for spin is {\it  vanishing of the second Stiefel-Whitney class} 
  $w_2(T)\in H^2(X;\mathbb Z_2)$; for instance, if $\pi_2(X)=0$, then the universal covering of $X$ is spin.} 
 It is also known that  all 3-manifolds are spin, while the complex projective plane  $\mathbb CP^2$ is non-spin. }  then $f$  {\it can't be strictly $dist^\natural$-decreasing}:
{\it there exists  a pair of points $x_1,x_2\in \partial X$,  such that 
$$dist^\natural(x_1,x_2) \leq dist^\natural(f(x_1),f(x_2)).
$$}}

{\it Remarks.} (a) As far as  the proof of  1.A is concerned, one needs only a very special case of this theorem, namely, where $X$ is  also a smooth convex domain in $\mathbb R^n$ and $Y$ is the unit ball in $\mathbb R^n$. Amazingly, however, even in this case, the only available  proof of 
2.A relies on the spin geometry and Dirac operators (see below).

(b) The assumption $Sc(X)\geq 0$ is, obviously, essential:  there is no curvature constrains on the boundaries of general Riemannian manifolds. 

But what is non-obvious, is how sensitive the geometry of $\partial X$ may be to the sign of the scalar curvature of $X$.

For instance, in agreement    
with   the {\it positive mass theorem} in general relativity,  there is no  Riemannian metric $g$  on the unit  ball $B^n\subset \mathbb R^n$ with $Sc(g)>0$  and 
with  $dist^\natural_g$ (non-strictly)  greater than the original $dist^\natural$ on the unit sphere  $S^{n-1}=\partial B^n\subset \mathbb R^n$, 

(c) It is unknown if the spin condition is essential.\vspace {1mm}

The second  components of the proof of 1.A - this is  an actual contribution of the present paper, is the following.

\textbf {2.B.   $dist^\natural$-Approximation Theorem.}\footnote{Smoothing the corners + "Dirac with boundary" is also used by  Simon Brendle    [Br2023]  in the proof  of  his polyhedral scalar curvature rigidity theorem.}  {\sf Let $X$ be a {\it compact  mean convex} Riemannian $n$-manifold with corners. Then, for all $\varepsilon>0$, there exists a smooth mean convex  hypersurface $\mathcal Y=\mathcal  Y_\varepsilon \subset X$ and a homeomorphism 
 $\psi =\psi_\varepsilon:\partial X\to \mathcal  Y$ with the following properties.

$\bullet_1$ The map $\psi$ is {\it $\varepsilon$-close to the identity, } $dist(\psi_\varepsilon(x),x)\leq \varepsilon$ for all $x\in \partial X$.

$\bullet_2$ The  $dist^\natural$ in $\cal Y$ is {\it greater than $dist^\natural_\rangle$ in $\partial X$ up to an  $\varepsilon$-error},
$$ dist^\natural_{\mathcal  Y}(\psi(x_1),\psi(x_2) )\geq (1-\varepsilon) 
dist^\natural_{\rangle}(x_1,x_2)$$
for all pairs of points  positioned within distances $\geq \varepsilon $  from the corners of $X$.}

We shall proof this in section 4, where we also show that, in the case of convex domains $X\subset \mathbb R^n$, the approximation
is possible with {\it strictly convex} $\mathcal Y$.

Then  we shall see  in section 5   that  2.A and 2.B (trivially) imply the following generlization of 1.A.
\vspace {1mm}

\textbf {2.C. Riemannian Angular Spread  Theorem.}  {\sf Let $X$ be a compact orientable {\it mean convex}  Riemannian $n$-manifold with  corners
and with {\it non-negative scalar curvature, $Sc(X)\geq 0$.} If  $X$ is {\it spin}, then the cubical $\rangle$-spread of  $X$ is universally bounded as follows,
$$\square^n_\rangle(X)\leq D=D_n\leq  2(n-1)\sqrt {n}.$$}

{\it Technical Strictness  Remark.}  Non-strictness of mean convexity may create      inconvenience, e.g. 
 a terminological one in dealing with vanishing "metrics".
But this is mainly  irrelevant,   since, in the cases of our immediate interest, e.g. for compact smooth hypersurfaces  
  in $\mathbb R^n$, strictness of mean convexity,   $mean.curv\geq 0\leadsto mean. curv>0 $, can is  easily achieved by arbitrarily   $C^\infty$-small perturbations. 

In general, with a minor analytic  effort,  one can $C^{\infty}$-approximate a   compact connected mean  boundary  
$ \partial X$ of a Riemannian manifold $X$ corners by a {\it strictly convex} hypersurfaces $Y\subset X$,
unless this  $\partial X$ is smooth  (no corners) with zero mean curvature.
 
 Thus,  one may    assume strictness of mean convexity in the present paper  whenever this helps  to  simplify understanding.

\vspace {1mm}


\subsection {$Sc$-Normalized Metric $g^\circ$,  Derivation  of 2.A from the LGSL Theorem  and $\rangle$-Capillary Problem}

 The 
counterpart of $g^\natural$  for Riemannian manifolds $X=(X,g)$ with positive scalar curvatures 
$$Sc(X,x)=Sc_g(x)>0$$ 
is the Sc-normalized  Riemannian metric 
$$g^\circ=g^\circ(x)=Sc_g(x)g(x)$$
on $X$.

 The basic geometric property of this  $g^\circ$ is the following special case of the 
Llarull -Goette-Semmelmann-Listing theorem (see 3.1.2 in [Gr2021] and references therein)\footnote{The essential ingredients  of the proofs in   [Ll1998], [GS2002], [Lis 2010] is  a sharp evaluation of eigenvalues of certain operators $\cal R$ in moduli over  Clifford algebras, where these $\cal R$ are algebraically associated with  the curvature operators $R$ of the underlying Riemannian manifolds $X$.  This suggests  a  direct Clifford algebraic approach to  the geometry of convex polytopes, where the complementary dihedral angle 
play the role of $R$ (compare with [WXY2022]).}

{\it \textbf {2.1.A. Euclidean  $dist^\circ$-Area Non-Contraction Theorem.}} {\sf Let $X$ be a connected  orientable  
$n$-dimensional Riemannian manifold     with $Sc(X) >0$ and let $\underline X\subset \mathbb R^{n+1}$ 
  a closed convex hypersurface. 
  
  Let $f:X\to \underline X$  be a smooth $g^\circ$-area decreasing map, that is 
  $$area_{\underline g^\circ} f(S)<area_{g^\circ} (S)$$
  for  all  smooth  surfaces   $S\subset X$, where $\underline g$ is the induced Riemannian metric in 
  $ X\subset \mathbb R^{n+1}$.}  

\hspace {2mm} {\it If $X$ is spin, then the map $f$ has degree zero  (hence, $f$ is contractible).}\vspace {1mm}

{\it Remarks} (a)  It is unknown,   not  even for $n=4$, if the spin condition is essential.

(b) The proof of 2.1.A depends on the index and vanishing theorems for the Dirac operator on $X$ with
 coefficients in  the vector bundle induced by $f$ from a unitary bundle on $\underline X$.
 
 The simplest kind of result of this kind, where the proof is technically very simple (see [GL1980], says that

 for {\it no
  Riemannian metric}  $g$ on $S^n$  the corresponding $g^\circ$  {\it can be significantly greater} than the 
   spherical metric:
 
\textbf {2.1.B.}  {\it If 
$dist_{g^\circ}\geq C dist_{S^n}$,
then $C\leq C_n$ for a universal constant $C_n$.}

(In fact, $C_n=\sqrt {n(n-1)}$  by Llarull's theorem [Ll1998].)

(c) If  $g$ has constant scalar curvature, then   2.1.B  (but not 2.1.A) can be proven by the technique 
of the geometric measure theory following ideas from  [SY1979].

Moreover: 

\textbf {2.1.C.}  {\it  If a metric $g$ on the unit ball $B^n\subset\mathbb R^n$ satisfies  $Sc(g)\geq C_n$, then 
the identity map $id :(B^n, g)\to (B^n, g_{Eucl})$
can't be distance decreasing.}

This is proven in [GL1983] for $n\leq 7$  and extended  to all $n$ in    [Gr2018] modulo [SY2017],   and directly   in [Loh2028].

(d) The obvious counterpart   of 2.1.A for open manifolds  fails to be true.
   
{ \textbf {2.1.D.  Example.}} The Euclidean space  $\mathbb R^n$, $n\geq 2$, admits a  Riemannian metric $g$  with $Sc(g)>1$
 and  such that $g^\circ$   is greater than the Euclidean metric.
$$g^\circ \geq g_{Eucl}.$$

(Notice that $g^\circ$  for such a $g$ is complete, but, (see [GL1983]),  $g$ can't be complete.)

{\it Proof.} Recall  that the scalar  curvature of the metric  $g_\phi=dx^2 +\phi^2 (x)dy^2 $
 on the $(x,y)$-plane is
  $$Sc( g_\phi(x))=-2 \frac {\phi''(x)}{\phi(x)}.$$
 Thus, if $\phi(x)$ is a strictly concave positive  function on the open interval $(0,1)$, such that   
 the integrals $\int_0^{1/2} \frac{\phi''(x)}{\phi(x) } dx$ and $\int_{1/2}^{1} \frac{\phi''(x)}{\phi(x)} dx$ diverge,
 then the metric $g^\circ_\phi$ on the band $U=(0,1)\times (-\infty\times\infty)$ is complete.
 Moreover, for all $
 \varepsilon>0$  there obviously exists a  distance decreasing   diffeomorphism from 
 $(U, g_\phi^\circ)$ onto $\mathbb R^2.$

Now  let $\phi(x)$ be equal  $x^\alpha$ near $x=0$  and to $(1-x)^\alpha$ near $x=1$  for $0<\alpha<1$,  
observe that these integral diverge and  make our example   with the obvious 
 distance decreasing  diffeomorphism $U\times \mathbb R^{n-2}\to  \mathbb R^{n}$.

{\it \textbf{On  Rediuction  of 2.A to 2.1.A}}.   This is achieved for  a manifold $X$ with a (mean convex) boundary by applying 2.1.A to the  
 double  \DD$(X)$ with a a suitably smoothed metric on it (see section 3.5 of [Gr2021])).
A more direct  but analytically more involved proof of 2.A  based on the  the index theorem for manifolds with boundaries was given in [Lott2020].

Then, on  the next    level of  sophistication,   the index theory    directly applies to  manifolds with corners  [WXY2022].

 This, formally speaking,   delivers   a two line proof of 1.A, but my  unsatisfactory  understanding of the techniques developed in [WXY2022] makes  me  reluctant to make such a shortcut in the proof.

\vspace {1mm}

{\it On  Capillary Geometry of $X$}. The above example highlights the difficulty   of applying the  geometric 
measure theory to $g^\circ$ and $g^\natural$, but it doesn't fully  rule out  such applications. 

Here is an instance of  what one may  expect of such an application.

 Let $X$ be a mean convex Riemannian $n$-manifold with corners and with positive scalar curvature  
and let $F^{n-1}_\mp \subset \partial X$ be two   faces positioned  {\it far away one from another in a suitable sense}, where the weakest  
condition (which may fail to be sufficient) would be a lower bound on the distance
 $dist^\natural_ \rangle $ between them:  $dist^\natural_\rangle(F^{n-1}_-, F^{n-1}_+)\geq const_n$, where,  ideally, 
 $const_ n= \pi$.

Then   one wants to have   a smooth hypersurface $Y\subset X$ with  $ \partial Y\subset \partial X$  transversal to the faces of $X$ 
and a smooth positive  function $\psi(y)$  on $Y$, such that the the  $\psi$-warped product  of $Y$ with the circle, 
$X_\rtimes =(Y\times \mathbb T^1, g_\rtimes)$,
for  $g_\rtimes =g_Y+\psi^2dt^2$, where $g_Y$ is the induced Riemannian metric in $Y$, 
 such that  the following conditions are satisfied:
 
$\bullet_{Sc}$  the metric $g_\rtimes$  has positive scalar curvature, 
 
$\bullet_{mean} $  the (boundary of the)  manifold $X_\rtimes $ is   mean convex,
 
$\bullet_{dist}$  the 
 $dist^\natural_\rangle$-distances between $(n-1)$-faces  in  $X_\rtimes$ are bounded from below,
  possibly, times  a controlled  $(1+\delta_n)$-factor, by the    $dist^\natural_\rangle$-distances between the corresponding
faces in $X$.

This would allow an inductive proof of (a sharp version?) of 1.A, where,
 observe,  
the expected $Y\subset X$,  say for $dim(X)= 3$  is a minimal surface (or something of this kind), 
which,  even
 for 3-polytopes $X\subset \mathbb R^3$ is by no means flat.
(Compare with [Gr2014'], [Gr2018], [Li2019] and section 5.81 in [Gr2021].)

\vspace {1mm}


\section {Rounding the Corners and  $dist^\natural$-Approximation of Simple    Polytopes}

Let $X\subset \mathbb R^n$ be a convex polytope and 
$\nu:\mathbb R^n \to X$ be the {\it normal projection}, that is   $\nu(x)\in X$  is  the nearest point point 
 to $X$, i.e. 
 $$dist  (x,\nu(x))=dist(x, X),  x\in\mathbb R^n,$$ 
and let 
$X_\circ=X_{\circ_\varepsilon} \supset X$, $\varepsilon>0$,  be  the $\varepsilon$-neighbourhood of $X$ that is    
 the set of points $x\in \mathbb R^n$  with $dist(x,X)\leq \varepsilon $.
 
Observe the following (compare with section 5.7 in [Gr2014] and 11.3 in [Gr2018]). 

 $\bullet_{\cup G_k}$  The boundary $\partial X_\circ$ is equal to  the union  of     {\it closures of the pullbacks of the open   $k$-faces} 
 $F^k\subset X$, $k=0,1,...,n-1$   intersected with $\partial X_\varepsilon$, denoted 
 $$ G_k=\nu^{-1}(F^k)  \cap  \partial X_\circ\subset \partial X_\varepsilon, $$
where such a $ G_k\subset \mathbb R^n= \mathbb R^k\times \mathbb R^{n-k}$ is  isometric  to the   product of the  corresponding face 
$F^k\subset \mathbb R^k $   by a  
convex $\varepsilon$-spherical polyhedron (dual to the normal section of  $ F^k $)  denoted  
$$ F_\perp^k \subset S^{n-k-1}(\varepsilon)\subset \mathbb R^{n-k}. $$ 

Thus,  {\it  the principal curvatures of} $G_k\subset \mathbb R^n$  are
  $$\underset {k}{ \underbrace {0,...0}},\underset {n-1-k}{\underbrace { \frac {1}{\varepsilon},..., \frac {1}{\varepsilon}}}\leqno\hspace{5mm} \bullet_{curv}$$
 and their mean curvatures 
     satisfy
   $$mean.curv(G_k))=\frac {1}{\varepsilon^{n-k-1}}.\leqno\hspace{5mm} \bullet_{mean}$$
 
 $\bullet_{C^1} $ 
Different 
 $G_k$, which  intersect  across parts of their boundaries,   have equal tangent spaces at their common points; thus the 
 boundary  
 $\partial X_\circ\subset \mathbb R ^n$  is $C^1$- actually $C^{1,1}$-smooth.

 {\it Quadratic Form  $g^\natural_{\circ_\varepsilon}$ and  Definition of $  dist^\natural_\circ=g^\natural_{\circ_\varepsilon}$.}  Let $g^\natural_{\circ_\varepsilon}$ be the product of the induced Riemannian metric 
 on the hypersurface  $\partial X_\varepsilon\subset \mathbb R^n$ by the squared  mean curvature of this hypersurface,
 $$g\natural_{\circ_\varepsilon}=(mean.curv)^2 g_{Eucl}|\partial X_{\circ_\varepsilon}$$
  and  observe that the metric defined with this Riemannian form $g^\natural_\varepsilon$ is exactly our 
  $dist^\natural_{\partial X_{\circ_\varepsilon}}$, which is denoted here $  dist^\natural_{\circ_\varepsilon}$. 

\textbf {3.A. $dist^\natural_{\circ}$-Convergence Theorem.} {\it If $X$ is a simple polytope, then the  $g^\natural_\varepsilon$-distance converges to the $\rangle$-distance,
$$dist^\natural_{\circ_\varepsilon}(K_1, K_2)\underset {\varepsilon \to 0}\to  dist_\rangle(F_1,F_2)$$
for all pairs of compact subsets in open $(n-1)$-faces $F_1,F_2\subset X$.}
$$K_1\subset F_1, K_2\subset F_2\subset X.$$

{\it $\square$-Example.} If $X=\square^n=[-1,1]^n$ is the $n$-cube, where, as we know, $dist_\rangle$ between 
opposite   $(n-1)$-faces 
is $ \pi$,  the  $g^\natural_{\circ_\varepsilon}$-distance between the corresponding faces in  $\square^n$ is only $\pi/2$.
To get the full $\pi$, one needs to go $\varepsilon$  away from the boundaries of these faces.

{\it Proof.}\footnote {Compare with section 5.7 in   [Gr2014] and 3.1.2(8) in [Gr2021]}. Let $Q\subset \mathbb R^n$ be a convex    polyhedral  $n$-dimensional cone  and $R\subset \partial Q_{\circ}=\partial Q_{\circ_1}$  be the complement to the flat part of  
$\partial Q_{\circ}$, that is the union of all $G_k$  with $k\neq {n-1}$. 

Observe that this  $R$  is a connected $(n-1)$-manifold with a boundary, where the connected components of this boundary are equal to the boundaries of the $(n-1)$-faces of $Q$.

\textbf {3.B. Minimal Path Lemma.} {\sf The shortest 
  paths $\gamma\subset R$  between  different connected components
 $\partial_1, \partial_2\subset  \partial R$  are    geodesic segments  contained in the  subsets 
 $G_k =F^k\times F^k_\perp\subset R\subset \partial Q_{\varepsilon=1}$, 
or in the intersection of several such subsets.}

Consequently, 

{\it the Riemannian  distance  between  $\partial_1$ and  $\partial_2$ is equal to the spherical  distance  between the intersection of 
 $\partial_1$  and $ \partial_2 $ with  the spherical polytope $G_0=Q_{\varepsilon=1}\cap S^{n-1}.$}
  $$dist(\partial_1, \partial_2)=dist_{S^{n-1}}(\partial_1\cap 	S^{n-1},\partial_1\cap S^{n-1}).$$

{\it Proof}. A priori,  $\gamma$ (which  is a $C^1$-smooth curve)  is composed  of several     geodesic  segments contained   in different $G_k$
But since {\it all} geodesic  segments in all $G_k$ 
 are {\it distance minimizing},  $\gamma$ is equal to the geodesic  continuation of  its  initial segment, say $\gamma_1$ 
 in some $G_k$; thus $\gamma$ stays in   this very $G_k$ all along. QED.

Now, let  a path $\gamma^\natural\in \partial X_\circ$ implement 
 the distance   $dist^\natural $  between two flat cells   in $ \partial X_\circ$, say between $G_{n-1,1}$ and $G_{n-1,2}$
The length of this path is equal to the sum of  $dist^\natural$ between components, say  $\partial_1$ and $\partial_1$,  of the boundary of the non-flat part
 $R\subset  \partial X_\circ$  crossed by $\gamma\natural$.
 
  If  $X$ is simple and  all $G_k=F_k\times \Delta^{n-1-k}$,  where  $\Delta^{n-1-k}$ are spherical simplices, these distances, because of 
  3.B, 
  can  implemented by geodesic segments  in $G_k$ with $k=n-2$ and 3.A  follows.

{\it \textbf { About Non-Simple  $X$ }}. Examples   show that     3.A fails to be true  for non-simple polyhedra $X$, but, due to  3.B  it allows  a a modification applicable to non-simple $X$.

Namely, the $(n-2)$-adjacency graph $E$ must be  replaced by the {\it full adjacency graph} $\mathcal E_+(X)\supset E(X)$, which, similarly  to $E$, has the $(n-1)$-faces  for vertices and  where the edges correspond to pairs of $(n-1)$-faces
which meet at $0$-faces (vertices) of $X$ and where the  lengths of these edges are defined by the corresponding    angles between these faces. 
 
However the resulting version of 1.A for non-simple $X$  doesn't bring anything new since 
it  follows from the "simple" case by a generic perturbation of the 
$(n-2)$-faces of $X$.  


\section {Locally Conical Hypersurfaces  and the $dist^\natural$-Approximation Theorem for  Non-simple $X$}

{ \it \textbf   {4.A. Conical Function Lemma.}} {\sf Let  $Y\subset \mathbb R^{m}$ be a (possibly unbounded, e.g. conical) convex polytope.  Then, for all $\varepsilon >0$  
there exists a {\it positive concave,} 
 function $ \phi=\phi_\varepsilon :Y\to \mathbb R_+$, which is piecewise  smooth   in the interior of  $Y$,  which {\it vanishes  on the boundary} $\partial Y$  and   which  satisfies the following four conditions.\footnote {With a little extra effort one can make $\phi$ smooth in the interior of $Y$.}  

$\bullet_\varepsilon$  The  {\it directional derivatives} of  $\phi$ at all {\it boundary points} $y\in \partial Y$ are {\it bounded} in absolute values by   $\varepsilon$, or equivalently 
$\phi$  is {\it $\varepsilon$-Lipschitz}: 
 $$|\phi(y_1)-\phi(y_2)|\leq\varepsilon \cdot dist(y_1,y_2)
\mbox  { for all }  y_1,y_2\in X.$$

$\bullet_{curv}$ The {\it principal  curvatures} of the graph  $\Gamma_\phi\subset \mathbb R^m\times \mathbb R_+$ at the smooth  points of $\phi$  are  everywhere {\it strictly positive.} 

$\bullet_{mean}$  The mean curvature of $\Gamma_\phi$ is uniformly positive on compact parts of $Y$ at  
 smooth points $(y,\phi(y))\in \Gamma_\phi$,
$$mean.curv( \Gamma_\phi, (y,\phi(y))\geq \epsilon(y)>0.$$
for a positive continuous function $\epsilon(y)$ on $Y$.

Moreover,  for all $(m-2)$-faces  $F\subset \partial Y$,
$$mean.curv(\Gamma_\phi (y, \phi(y)))\geq const\frac {1}{dist(y, F)}\leqno{ \hspace {6mm}\bullet_{1/d}}$$
for some positive constant $const=const_{P, \varepsilon}>0$ and all interior points $y\in Y$, where $\phi$ is smooth.}

{\it Proof.} The existence of $\phi$ is obvious for $m=1$ and the general case follows by induction in $m$ as follows. 
 
Represent $Y$ by the intersection of the wedges $W_i\subset \mathbb R^m$, $i=1,...,j$, 
which are based on the  1-faces  $F_i^1\subset Y$,
$$ Y=\bigcap_{i=1}^jW_i, \mbox { for }   W_i=C_i\times L_i,\mbox    { and $C_i\subset \mathbb  R^{m-1}_i$},$$
where 

$\bullet$  $L_i\subset \mathbb R^m$ are the  straight  lines, which extend the 1-faces $F_i^1\subset Y$;
 
 $\bullet$ $\mathbb R^{m-1}_i\subset \mathbb R^{m}$  are normal spaces to the faces $ F_i^1$ at some points $y_i\subset F_i^1$;
 
  $\bullet$ $C_i\subset \mathbb R^{m-1}_i$  are the convex  tangent cones to $Y$ at the points $y_i$,  that are the conical extensions of the intersections of 
 $ \mathbb R^{m-1}_i$  with small neighbourhoods of $y_i$ in $Y$.

Let $\phi_i(c)$ be concave   functions in the cones $C_i$, which satisfy the four conditions  $\bullet_\varepsilon$, $\bullet_{curv}$, 
$\bullet_{mean}$,  $\bullet_{1/d}$, 

 let $\bar\phi_i(c, l)=\phi_i(c)$
for $(c,l)\in W_i=C_i\times L_i\subset \mathbb R^m$

and let $\bar\phi$ be the minimum of generic $\lambda_i$-perturbations  of the functions $\bar\phi_i$  on $Y$,
$$\bar\phi(y) =\min_i \lambda_i \phi_i(y), y\in Y.$$
for small  {\it generic}  $\lambda_i>0$. \footnote {Generic  $\lambda_i$ are needed to assure    piecewise smoothness
 of $\bar\phi(y)$ in the interior of $Y$.} 

Clearly, the function  $\bar\phi$ satisfies  $\bullet_\varepsilon$, $\bullet_{curv}$, 
$\bullet_{mean}$, but it may fail   $\bullet_{1/d}$  at the vertices  $y_\nu\in Y$. 

To correct this,  modify $\bar\phi$ at  $y_\nu$ as follows.
Let $U_\nu\subset Y$ be (very)  small (pyramidal)  neighbourhoods of $y_i\in Y$, which are  bounded in $Y$ by hyperplanes cutting $y_\nu$ away  from $Y$,
let 
$$\bar Y=Y\setminus \bigcup_\nu U_\nu$$
be the correspondingly truncated $Y$ and let $\phi(y)$ be the {\it smallest concave} function on $Y$,
which is equal to  $\bar \phi$ on $\bar Y$ and which vanishes on the boundary of $Y$.

In geometric term, the convex body  $Y^+_\phi \subset \mathbb R^m\times \mathbb R_+$  under the graph   $\Gamma_\phi \subset \mathbb R^m\times \mathbb R_+$ is obtained by firstly cutting away $y_\nu$ from  $Y^+_{\bar \phi} \subset \mathbb R^m\times \mathbb R_+$   by vertical half-hyperplanes $H^+_\nu\subset  \mathbb R^m\times \mathbb R_+$
and then adding the cones from $y_\nu$ over the intersection $Y^+_{\bar \phi}\cap H^+_\nu$ to the resulted truncated $Y^+_{\bar \phi} $ 
.

Now, clearly,  the mean  curvature of $\Gamma_\phi$ does blow-up  as $1/d$ for the distance $d$  to the $(m-2)$-faces of $Y$ and the proof of 4.A is concluded.
\vspace {1mm}

{\it \textbf  {Proof of the  $dist^\natural$-approximation theorem  2.B for convex polytopes.} }  Let $X\subset \mathbb R^m$ be a compact 
convex polytope  and  let $X^+_\varepsilon\subset \mathbb R^n$ be obtained by adding the subgraphs of
 the functions  $\phi =\phi_\varepsilon$ on all $
(n-1)$-faces $Y$  of $X$   to $X$.
 
 The following five properties of $X^+_\varepsilon $ trivially follow from 4.A.

$\bullet_\delta$ The set $X^+_\varepsilon$ is {\it pinched} between $X$ and a (small) $\delta$-neighbourhood of $X$,  
$$X\subset X^+_\varepsilon  \subset U_\delta(X) \subset  \mathbb R^n,\mbox {  where $\delta\to 0$ for $\varepsilon \to 0$}.$$

$\bullet_{conv}$ If $\varepsilon>0$ is sufficiently small, then  $X^+_\varepsilon$ is {\it convex}.

$\bullet_{n-2}$  The intersection of the boundary of   $X^+$ with  $X$ is equal to the {\it union of the closed $(n-2)$-faces} of $X$, 
$$  \partial X^+_\varepsilon \cap X=  \partial X^+_\varepsilon\cap \partial X=\bigcup_{i=1,..., k\leq n-2} F^k_i.$$

$\bullet_\angle$ The {\it dihedral angles} of  $X^+_\varepsilon$  along $(n-2)$-faces of $X$ (contained in   $\partial X^+_\varepsilon$)\footnote{These are  the  angles between the pairs of extremal supporting hyperplanes  to  $X^+_\varepsilon$ at the points $x\in F^{n-2}\subset X\cap X^+_\varepsilon$.}   are   bounded by the dihedral angles of $X$ between these faces as follows,
$$\angle^+\leq \angle+2\varepsilon.$$

$\bullet_{1/d}$ The mean curvature of $\partial X^+_\varepsilon$ at smooth points 
$x\in  \partial  X^+_\varepsilon$  satisfies
$$  mean.curv   (\partial  X^+_\varepsilon, x)\geq const\frac {1}{dist(x, F^{n-3})}$$
for some $const>0$  and all $(n-3)$-faces $F^{n-3}$ of $X$.

It follows, that paths $\gamma\subset \partial X^+_\varepsilon $,  which approach  $F^{n-3} $ have {\it infinite} 
$g^\natural$-lengths; hence $g^\natural$-shortest paths cross $(n-2)$-faces away from $(n-3)$ faces.  

Then an additionally $C^2$-smoothed  boundary $Y$ of 
the $\varepsilon_\circ$-neighbourhood  $U=U_{\varepsilon_\circ}(X^+_\varepsilon )$ serves as the required approximation of $X$
by the (trivial)  argument  from 3.A.

{\it \textbf{ Generalization to  Mean Convex  Manifolds $X$ with Corners.}} Think of $X$ as a mean convex domain with corners 
 in a larger Riemannian manifold, say   $W\supset X $ and construct $X^+_\varepsilon\subset W$ 
 in three steps. 
 
 1. Make the  $(n-1)$-faces $F_i^{n-1}$ of $X$ {\it strictly} mean convex by  $C^\infty$-perturbations, 
 while keeping these faces unperturbed   on  the parts of their boundaries which are  close to $(n-3)$-faces, i.e. on the intersections 
 $\partial F_i^{n-1}\cap  U_\varepsilon (F^{n-3}_j)$.
 
 This is done by linearizing the problem as it is done in the first proof of ({\Large $\star$}$_>$) in  section 11.2 of 
 [Gr2018]\footnote{This argument is outlined in   [Gr2018] for {\it simple}  $X$, i.e. where the $(n-1)$-faces intersect
  transversally  
and thus   the combinatorial structure of $X$ remains stable  under small perturbations, while in the present case, one needs to keep
the perturbation fixed on $\partial F_i^{n-1}$ near the $(n-3)$-faces to preserve the combinatorial structure of $X$. 
 
 In any case, all this is a minor matter and one doesn't loose  much by assuming that $X$ is strictly mean convex to start with.

Also notice  that second "variational proof"  of     ({\Large $\star$}$_>$)  in   [Gr2018] is invalid.}

 {\it Warning.}   One can't, in general,  achieve this while keeping the  faces  fixed everywhere  on their boundaries as it was done for convex $X\subset \mathbb R^n$.
  
For instance, if $X$  is a locally convex geodesic  polygon in a Riemannian  surface $W$,
then  an  edge $F^1$ in $X$ can be approximated by a strictly convex curve  with the same ends as $F^1$,
if and only if  $F^1$,  which   itself  is a geodesic segment, contains {\it no conjugate points}.
 
2.  At the  second step one make the mean curvature of the faces blow up at the $(n-3)$-faces  with the rate $1/d$.
   as in the above convex case. In fact, since this blow-up  property is invariant under diffeomorphisms, one can   perform it locally 
   in normal geodesic coordinates and then  glue these together by  a partition of unity argument. 

Notice that this is unneeded if $X$ is simple, where one goes directly to the third step.

3. Once 1 and 2 are done and one arrives at a   {\it strictly} mean convex $X_\varepsilon $,
   which satisfy the above $\bullet_{1.d}$, then,
as earlier, one takes the  $C^\infty$-smoothed boundary of a  small  $\varepsilon_\circ$-neighbourgood  
$U_{\varepsilon_\circ}(X_\varepsilon) $ for the required approximation $\mathcal Y$ of $\partial X$  
(compare with 5.7 in [Gr2014] and 11.4 in [Gr2018]).

{\it Remark.} It would take a couple of extra  pages to explicitly  write down the (quite boring) details of the above
argument but it would add nothing new to what we  have already seen  in the convex case.
ing 

{\it \textbf {Convexly Stratified Manifolds.} } The step 2 in the above  argument  takes  $X$ out of 
the  category of 
manifolds with corners, where the new  manifolds are locally diffeomorphic not to convex polytopes 
but  to certain 
 smoothly stratified   convex subsets 
 $\underline X \subset \mathbb R^n$, such, for instance,  as cones over smooth convex bodies 
  in $\mathbb R^{n-1} \subset \mathbb R^n.$

The most general class $\underline { \mathcal X}_{gen} $ of such  $\underline X$,  where  the statement of theorem 2.A
makes sense, consists   of closed convex  domains  $\underline X$, such that      the boundaries  of $X$  
are  piecewise smooth  in the complements of   closed $(n-3)$-dimensional subset $Z\subset \partial \underline X$.

Probably, the proof of  theorem 2.A can be extended to the corresponding class $\mathcal  X$  of mean convex Riemannian
manifolds $X$ locally diffeomorphic to such $\underline X$. 




\section{Lipschitz Maps and the Proof of Theorems 1.A. and 2.C.}

Here is an essential, albeit    elementary (and trivial), geometric  fact  one needs. 

\textbf {5.A. Lipschitz Mapping Lemma.}  Let $Y$ be a closed orientable Riemannian $(n-1)$-manifold and $\phi$ be a continuous map
from $Y$ to the boundary of the $n$-cube $\square^n= [-1,1]^n$,
$$\phi:Y\to \partial \square^n,$$ 
such that the distances between the  pullbacks of the opposite faces $\square^{n-1}_{i\pm}\subset \square^n$, $i=1,...,n$,
satis
$$dist (\phi^{-1}(\square^{n-1}_{i+}, \phi^{-1}(\square^{n-1}_{i-})) > D.$$

Then the composition of $\phi$ with the obvious radial homeomorphism from  $\partial \square^n$ to the unit sphere
$$ \partial \square^n\to S^{n-1} \subset \square^n$$
is {\it homotopic to  a smooth map,
$$f : Y\to S^n,$$  
such that the differential of $f $  satisfies 
$$||df||<\frac {2 \sqrt {n}}{D}.$$}

{\it Proof.} Let $\delta_i(y)$ be the distance functions to   $\phi^{-1}(\square^{n-1}_{i-})\subset  Y$ truncated by $D'>D$ such that
$$D<D'< \max_i dist (\phi^{-1}(\square^{n-1}_{i+}, \square^{n-1}_{i-}),$$ 
namely, 
$$\delta_i(y)=\min(dist (y, \phi^{-1}(\square^{n-1}_{i-})), D')$$
and observe that the map 
$$\Delta=\left (\frac{2}{D'} \delta_i(y)-1,..., \frac{2}{D'}\delta_n(y)-1\right)$$
sends 
$$Y\to \partial \square^n\mbox {  for  }\square^n=[-1,1]^n \subset \mathbb R^n,$$ 
that  this map is {\it homotopic to} $\phi$ 
and that it is
 $\frac{2\sqrt{ n} }{D'}$-{\it Lipschitz. }
Since the radial map radial map  $\partial \square^n\to S^{n-1}\to S^{n-1}$ is  distance decreasing  and $D'>D$,
the composed map 
$$Y\overset {\Delta}\to \square^n \to S^n$$
can be approximated by the required $f$.

{\it \textbf {5.B. Conclusion of the Proof of Theorem  2.C.} } Let $X$ be a Riemannian  manifold with corners as in 2.C, let   $\square^n_\rangle(X)>D $ and let
$\Phi:X\to \square^n$ be a  continuous combinatorial   map, which  
 satisfies  $\bullet_{dist_\rangle}$  and $\bullet_{deg}$ from section 1 and also  $\bullet_{dist_\rangle}$, but now  with $D$  instead of $d$.
 
  Then   a smooth   mean convex hypersurface  $\mathcal Y
  \subset X$, which approximates $\partial X$  according to  2.B,  which we assume strictly mean convex and which we  endow with the metric $dist^\natural$, and the map $\phi=\Phi_{|\mathcal Y}:\mathcal Y\to \partial \square^n$
  satisfy the assumptions of 5.A.  
  Hence, 
 $(\mathcal Y, dist_\natural) $  admits a   $\lambda$-{\it Lipschitz} map to the unit sphere 
$S^{n-1} $ for $\lambda<\frac{2\sqrt{ n} }{D}$  as in 5.A.

Since the degree of this map doesn't vanish according to $\bullet_{deg}$, theorem 
  2.A   says that   $\lambda\geq \frac {1} {mean.curv(S^{n-1})}=\frac {1}{n-1}$, 
  which implies that $\frac{2\sqrt{ n} }{D}\leq \frac {1}{n-1}$
   and $$D<2(n-1)\sqrt {n}.$$
QED


\section {Combinatorial Waists and $\rangle^{n-k-1}$-Angles}

 {\it The $F^k$-overlaps} of a  map  from a manifold $X$  with corners, e.g. from a polytope, to some set, 
 say $\alpha:X\to \Xi$,   denoted
$$\mbox {${\overset {\smile}\#}$$_\alpha ^k (X)$  and $ \widehat{\#}_\alpha ^k (X)$,  }$$
are  the maxima of the  numbers  of open, respectively closed, $l$-faces in $X$ the $\alpha$-images of which in  $\Xi$ have a common point $\xi$,

For instance, generic linear maps  $\alpha$  from 
$n$-polytopes 
$X\subset \mathbb R^n$ to 
$\mathbb R^{n-1}$
satisfy  
 {${\overset {\smile}\#}$$_\alpha ^{n-1} (X)=2$} and, 
if $X$ simple, then $  \widehat{\#}_\alpha ^{n-1} (X)=n+1.$

{\it \textbf {The spherical $(n-k-1)$-(co)angle} }  of  a  convex  subset  $X\subset \mathbb R^{n}$, e.g. a Euclidean $n$-polytope, at  a point $x\in\partial X,$
 denoted 
 $$\rangle^{n-  k-1}_x(X),$$  
is the $(n-k-1)$-dimensional spherical volume (Hausdorff measure) of the set of the supporting 
planes to $X$ at $x$, and  where we denote
$$\rangle_{F^{k}} =\rangle^{n-k-1}_{x\in F_\circ^k} (X)$$
 for open and closed  $k$-faces $F^{k}$ and  for   
 $F_\circ ^k\subset F^k$ being  the interior parts of these faces.

 Then the  "angle"  $ \rangle^{n-k-1}_x(X)$ at a point $x$ in a Riemannian   manifold $X$ with corners is defined as the corresponding   angle of the   tangent cone of $X$ at $x$

For instance, if $x\in F^{n-2}$,   this is the complementary dihedral angle of the face $F^{n-2}$  defined earlier.

Define ${\overset {\smile}\sum}\rangle^{n-k-1}_\alpha (X)$  and  $\widehat \sum\rangle_\alpha^{n-k-1} (X)$ for Riemannian manifolds $X$ with corners 
as the supremum over $\xi\in \Xi$ of the sums of these angles over the set of non-empy 
intersections of the $k$-faces $F^{k}$  in $X$ with the  
$\alpha$-pullbacks of points $\xi$ in $\Xi$,
$$\sup_{\xi\in \Xi} \sum  \rangle_{F^k},  \mbox  { $ F^{k}\cap \alpha^{-1}(\xi)\neq \emptyset$,}$$
i.e. the sum is taken over all  open, respectively closed, $k$-faces $F^{k}\subset X$,which intersect  $\alpha^{-1}(\xi)$.

{\it Remark.} This definition makes sense for all  weight functions $w$ on the faces instead of $\rangle^{n-k-1}$, where, e.g. for $w(x)=1$,   
one recaptures the numbers $\overset {\smile} \#$$_\alpha^{k} (X)$ and $\widehat {\#}^{k}_\alpha (X).  $

{\sc Problem.} Given a class $\mathcal A$ of function $\alpha$,  evaluate possible values  $ \#^k_\alpha (X)$ and 
$\sum\rangle_\alpha^{n-k-1} (X)$   for convex polytopes and other "interesting"
 manifolds with corners in terms of other geometric invariants.

{\it  \textbf  {Example  6.A.}} Let  $X\subset \mathbb R^n $ be a convex polytope and  $\alpha:X\to \mathbb R^{n-2}$  a {\it continuous} map.
If $X$ is simple,\footnote{This is probbaly redundant} then 
the number $\widehat \#^{n-2}_\alpha (X)$ is bounded from below by the combinatorial $\square$-spread of $X$ as follows
$$\widehat \#_\alpha^{n-2} (X)\geq cost_n\cdot \square^n_{comb}(X).$$

{\it Sketch of the Proof.}  Let  $g_\varepsilon$ be a   Riemannian metric on  $\partial X$, which distance-wise     $\varepsilon$-approximates    $dist_{comb}$ on $X$.
By the argument from the previous section, $(\partial X, g_\varepsilon)$  admits a 1-Lipschitz map $\Phi$   of non-zero
 degree to the sphere $S^{n-1}(R)$
of radius $R\geq const'_n   \square_{comb}(X)$. It follows by the (quite elementary)  1-waist inequality for spheres (see [Guth2014] and references therein) the $\Phi$-image of the pullback $\alpha^{-1}(\xi)$, $\xi\in \mathbb R^{n-2}$, has length $\geq 2\pi R$.
Hence, the $g_\varepsilon$-length  of $\alpha_{-1}(\xi)$ is also  $\geq 2\pi R$,  which,  since $X$ is simple,  implies the required
bound $\widehat \#_\alpha^{n-2} (X)\geq const'_n R$ for $\varepsilon\to 0$.

Probbaly, a similar argument applies to continous maps $\alpha:X\to \mathbb R^k$ for all $k=1,...,n-2$,
 thus showing,  at least for manifolds $X$  with simple corners,  that
$$\widehat \#_\alpha^k (X)\geq cost_n\cdot \left (\square^n_{comb}(X)\right)^{n-k-1}$$  
for all continuous  $\alpha$. 

But it is unclear what happens to  $\sum\rangle_\alpha^k$. \vspace{1mm}

\hspace {8mm}{\sf For instance, let $X\subset \mathbb  R^n$  be  a convex polytope and  $k=1,...,n-2$. }

{\it \textbf  {Question 6.B.}} Is 
$$\inf_\alpha{\overset {\smile}\sum}\rangle_\alpha^k (X)\leq cost_n$$ 
the infimum is taken over all continuous (may be even linear?)
maps $\alpha: X\to \mathbb R^k$?
\vspace{1mm}

{\it \textbf  {Question 6.C.}}  Does there exist an  $(n-k)$-dimensinal affine subspace $A\subset \mathbb R^n$, which transversally meets 
$N>0$  (open) $k$-dimensional  faces    $F_i^{k}\subset X$, $i=1,...N$, such that 
$$\frac {1}{N}\sum_{i=1}^N\rangle^{n-k-1} _{F_i^{k}}(X) \leq \frac{const_n}{(\square_{comb}(X))^{n-k-1}}?$$

The positive answers to this would yield  the  following  generalization of corollary 1.B to $k\leq n-3$.

{\it \textbf  {Conjecture  6.D.}} {\sf If the combinatorial $\square^n$-spread of a convex polytope  $  X\subset \mathbb R^n$ is large, then 
there exists a $k$-dimensional face $F_{min}^{k}\subset X$ 
with small  $\rangle^{n-k-1}$-angle:
$$\rangle^{n-k-1}_{F^k_{min}}(X)\leq const_n \left (\square^n _{comb}(X)\right)^{-(n-k-1)}, $$
or,  at least, 
$$\rangle^{n-k-1} _{F^k_{min}}(X) \to 0 \mbox { for }\square^n _{comb}(X)\to \infty.$$}
for simple polytopes $X$.

\vspace {1mm}


\section{ Surgery with Corners and Related Problems}

It is claimed in section 1.3   of [Gr2014']  that the so called {\it staircase  thin  surgery}  of mean convex manifolds 
with positive scalar curvatures   \footnote {See 
[GL1980'],  [BaDoSo2018], [Gr2021].       }can be also applied to manifolds $X$ with corners. 
However, I overlooked the difficulty in proving the following.

   {\it \textbf {7.A.  $\angle$-Shrinking Problem.}}  Let $Y_0\subset S^{n-1}\subset \mathbb R^n$ be a convex spherical polytope.
 
 {\sf  Does there exist a continuous deformation 
      $Y_t\subset S^{n-1}$, $0\leq t\leq1$, of  $Y_0$,  where all  $Y_t$ for $t<1$ are convex spherical  polytopes combinatorially  isomorphic  to 
      $X_0$ and having their  
       dihedral angles bounded by the corresponding angles   of $Y_0$  and where $Y_1$ is a single point? }   
   
 It is easy to construct    such a $Y_t$  for  $dim(Y_0)=2$, and also for "sufficiently round"   spherical polytopes of dimensions>2,  
      where such shrinking can be achieved by  projective transformations   of $Y_0$,   but I was unable to prove or disprove it for general $Y\subset  S^{n-1}$ if $n\geq 4$.

  And  granted such a deformation  for  the spherical base  $Y_0$ of the tangent cone  $TC_{x_0}(X)\subset  T_{x_0}(X)=\mathbb R^n$, 
at a vertex $x_0\subset X$,   say for a strictly  mean convex domain $X\subset \mathbb R^n$ with corners, 
 the  staircase construction  delivers another strictly  mean convex $X_1\subset \mathbb R^n$, such that 
  
  $\bullet_{cut}$ the domain $X_1$ is diffeomorphic to $X$ with the vertex $x_0$ cut away by a hyperplane 
  $H_0\subset \mathbb R^n$ parallel  to a supporting hyperplane of $X$ at at $x_0$, where  this diffeomorphism moves all points at most by a given $\varepsilon>0$  and 
    fixes the points $\varepsilon$-far from $x_0$;  

   $\bullet_{>\rangle}$ the  dihedral angles at the (old)  $(n-1)$-faces of $X_1$ away from the  cut $X\cap H_0$ 
 are bounded by the corresponding angles of $X$;

   $\bullet_{\pi/2}$  the   dihedral angles between  the new $(n-1)$-face corresponding to $X\cap H_0$  with the  old ones    are equal to $\pi/2$.

  Observe that if $n=3, $ this construction, when applied to 
all vertices of $X$, delivers a simple polytope and thus provides an alternative  reduction of the general case of theorem 1.A 
to that for simple $X$.  

However,  since 7.A.  remains problematic for $dim(Y)\geq 3$  

\hspace {10mm} {\it the thin surgery at the corners remains problematic as well. }

Also  pondering over 7.A brings to one's mind the following more  general problems.

  {\it \textbf {7.B.  $\rangle$-Variation Problem.}} Find the homotopy type   of the
   space $\mathcal X (\mathcal C,\kappa)$   of (possible)  dihedral angles  
  of convex $n$-polytopes $X$  of  given combinatorial type $\mathcal C$  in the space of constant curvature $\kappa$ 
  and determine   how this space varies depending on $\infty <\kappa<\infty$.

{\it \textbf {7.C. Scalar Curvature  $\rangle$-Problem.}} Let $X$ be a compact connected  smooth manifold with corners, let 
$-\infty <\mu_i<\infty $  
be numbers associated to the  $(n-1)$-faces of $X$ and $0<\alpha_j <\pi$   be associated to the $(n-2)$-faces.  
Determine the homotopy type of the space $\mathcal G(X, \sigma, \mu_i, \alpha_j)$, $\sigma>0$  of  Riemannian metrics $g$  on $X$ 
such that 

$\bullet_\sigma$  the  scalar curvature of $X$ satisfies: 
$$Sc(X)> \sigma; $$

$\bullet_\mu$ the mean curvatures  of the  $(n-1)$-faces of $X$ satisfy:
 $$mean.curv_g(F_i^{n-1}) >\mu_i;$$

$\bullet_\alpha$  the complementary   dihedral angles  at the $(n-2)$-faces satisfy: 
$$\rangle_g((F^{n-2}_j)\geq \alpha_j.$$

Also determine how this space varies depending on $(\sigma, \mu_i, \alpha_j)$.


\section { On Random Polytopes}

Let $\Sigma= \{\sigma_i\}_{i=1,...,N^{n-1}}\subset S^{n-1}  $ be randomly chosen points on the unit  sphere and 
$X_N=X(\Sigma)$ be the (necessarily simple)  convex polyhedron defined by the tangent hyperplanes to the sphere at the points $\sigma_i$. 

Let  $dist_{comb, N}(s_1,s_2)$, $s_1,s_2\in S^{n-1}$  be the   combinatorial distance between the $(n-1)$-faces $F_1,F_2\subset X$ 
of $X_N$ the normal projections of which to $S^{n-1}$ contains the points $s_1$ and $s_2$ respectively. (Never mind the distinction between open and closed faces.) 

\textbf {8.A. Spherical  $dist_{comb}$-Conjecture}.  {\sf There exists a universal constant $\Delta_n$ such that 
$$ \frac {dist_{comb, N}(s_1,s_2)}{N\cdot dist_{S^{n-1}}}\to \Delta_n\mbox  { for }  N\to \infty$$  
with probability 1 for all pairs of points $s_1,s_2\in S^{n-1}$:

the probability of the inequality $\left| \frac {dist_{comb, N}(s_1,s_2)}{N\cdot dist_{S^{n-1}}}-\Delta_n\right|>\varepsilon$ tends to zero for $ N\to \infty$  for all   $\varepsilon >0$.}

{\it Remark.}  Probbaly, this follows by the  results/arguments  from  [BDGHL2021] but I haven't looked at this closely.
\footnote {There is an extensive literature on random polytopes, where  much of  known estimates of the sizes of random  polytopes  concern {\it upper} bounds on  combinatorial edge-diameters, which are   motivated by the {\it Hirsch conjecture,} while we are  interested on lower  bounds on the $\square$-spreads.}
in any case 
 an elementary (Poisson)  percolation  argument  shows   that,  with  overwhelming probability,
$$\mbox {$ \frac {dist_{comb, N}(s_1,s_2)}{N\cdot dist_{S^{n-1}}}\leq cost_n$   and 
 $ \frac {dist_{comb, N}(s_1,s_2)}{N\cdot dist_{S^{n-1}}}\geq \frac {cost'_n}{\log N}$   for $N\to \infty.$} $$

\textbf {8.A. Spherical $dist_\rangle$-Conjecture}. {\sf  Let  $dist_\rangle(s_1,s_2)=dist_\rangle (F_1,F_2)$ for the above $F_1,F_2$. Then   
there exists a universal constant $\Delta^\rangle_n$ such that 
$$ \frac {dist_{\rangle}(s_1,s_2)}{dist_{S^{n-1}}}\to \Delta^\rangle_n\mbox  { for }  N\to \infty$$  
with probability 1 for all pairs of points $s_1,s_2\in S^{n-1}$.}

{\it Remark.}  Exact evaluation of $\Delta _n$ and   $\Delta^\rangle_n$ may be difficult but the ratio  $\Delta _n/\Delta^\rangle_n$
seems computable.

There  are other commonly used   definition of  "random polytope"  (see [Schneider2008]);    we single out the following. 

Let  $\mathcal C(n,M)$ be the set of combinatorial types of simple $n$-polyhedra $X$ with $M$ faces and observe that
the cardinality of this set is pinched between two exponentials:
$$A^M \leq \#\mathcal C(n,M)\leq B^M$$
Cutting $X$ by hyperplanes in two parts suggests that   $\log \#\mathcal C(n,M)$ is (essentially) super-additive,
and the limit 
$$\lim_{M\to \infty}\frac {\log \#\mathcal C(n,M)}{M},$$ 
which seems an interesting number, exists.

Then we assign equal probabilities to all points (combinatorial types) in  $\mathcal C(n,M)$ and conjecture that 
 
{\sf the graphs $E=E(X)$ of the   so defined random $n$-polytopes $X$ with $M$ faces endowed  with metrics  
${M^\frac {-1}{n-1}}dist_{comb}$ {\it Hausdorff converge} to the sphere $S^{n-1}(R_n)$ of some radius $R_n$  for $M\to\infty$.}

\section {References}

\hspace {4mm} [AK]   Karim Adiprasito, {\sl Private Communication.} \vspace{2mm}

 [BaDoSo2018] J. Basilio ,J. Dodziuk, C. Sormani, {\sl  Sewing Riemannian Manifolds with Positive Scalar Curvature,} The Journal of Geometric 
AnalysisDecember 2018, Volume 28, Issue 4, pp 3553-3602.

\vspace{2mm}

[BDGGL2021]  Gilles Bonnet,
Daniel Dadush,
Uri Grupel
Sophie Huiber,
Galyna Livshyt, {\sl Asymptotic Bounds on the Combinatorial Diameter of
Random Polytopes}  	arXiv:2112.13027.

\vspace {2mm}

[Br2023] S. Brendle,
  {\sl Scalar curvature rigidity of convex polytopes}
arXiv:2301.05087 

\vspace {2mm}




[GL1980] M.Gromov, B Lawson, {\sl Spin and Scalar Curvature in the Presence of a Fundamental Group I}, Annals of Mathematics, 111 (1980), 209- 230.
\vspace {2mm}

  [GL1980'] M.Gromov, B Lawson,     {\sl The classification of simply connected manifolds of positive scalar curvature},  Ann. of Math.  111 (1980) pp. 423-434. 
\vspace {2mm}

[GL1983] M.Gromov, B Lawson, {\sl Positive scalar curvature and the Dirac operator on 
complete Riemannian manifolds}, Inst. Hautes Etudes Sci. Publ. Math.58 (1983), 83-196.

\vspace {2mm}

[Gr2014] M. Gromov,{\sl  Plateau-Stein manifolds,} Central European Journal 
of Mathematics, Volume 12, Issue 7, pp 923-95.
\vspace {2mm}

[Gr2014] M. Gromov, {\sl   Dirac and Plateau billiards in domains with corners}, Central European Journal of Mathematics, Volume 12, Issue 8, 2014, pp 1109-1156.\vspace {2mm}

[Gr2018]  M. Gromov,  {\sl Metric Inequalities with Scalar Curvature,} Geometric and Functional Analysis Volume 28, Issue 3, pp 645-726.
\vspace {2mm}

[Gr 2021] M. Gromov, {\sl  Four Lectures on Scalar Curvature}arXiv:1908.10612 [math.DG].

\vspace {2mm}

[GS2002] S. Goette and U. Semmelmann, {\sl Scalar curvature estimates for compact symmetric spaces}. Differential Geom. Appl. 16(1):65-78, 2002.

\vspace {2mm}

[Guth2014] L Guth, {\sl The waist inequality in Gromov’s work}, MIT Mathematics math.mit.edu/~lguth/Exposition/waist.pdf

\vspace {2mm}

[Li2019] Chao Li, {\sl  
The dihedral rigidity conjecture for n-prisms}. https: // arxiv. org/ abs/ 1907. 03855 , 2019
\vspace {2mm}

[Lis2010] M. Listing, {\sl Scalar curvature on compact symmetric spaces}. arXiv:1007.1832, 2010.

\vspace {2mm}

[Lla1998] M. Llarull, {\sl  Sharp estimates and the Dirac operator}, Mathematische Annalen January 1998, Volume 310, Issue 1, pp 55- 71.

\vspace {2mm}

[Loh2018] J. Lohkamp {\sl  Contracting Maps and Scalar Curvature}
 arXiv:1812.11839

\vspace {2mm}

[Lott2020] J. Lott, {\sl Index theory for scalar curvature on manifolds with boundary,} arXiv:2009.07256

\vspace {2mm}

[MSE] Mathematics Stack Exchange {\sl How do dihedral angles grow with number of edges in Euclidean polyhedra}.

\url {https://math.stackexchange.com/questions/3732790/} \vspace {2mm}

[Sh2008] Rolf Schneider, {\sl Recent Results on Random Polytopes}
Bollettino dell'Unione Matematica Italiana (2008)
Volume: 1, Issue: 1, page 17-39

\vspace {2mm}

[SY(structure) 1979] R. Schoen and S. T. Yau,  {\sl  On the structure of manifolds 
with positive scalar curvature,} Manuscripta Math. 28 (1979), 159-183.\vspace {2mm}

[SY(singularities) 2017] R. Schoen and S. T. Yau
 {\sl Positive Scalar Curvature and Minimal Hypersurface Singularities}. arXiv:1704.05490 \vspace {2mm}

 [WXY2022] Jinmin Wang, Zhizhang Xie, Guoliang Yu,
{\sl On Gromov's dihedral extremality and rigidity \vspace {2mm}
conjectures}  arXiv:2112.01510v4 [math.DG]

\end {document}